\documentclass{article}
\usepackage{amsmath,amsthm}
\usepackage{graphicx,epstopdf,epsfig,multirow,epic,bm}

\normalsize \rm
\begin{document}

\begin{center}
{\Large  \textbf { Determining Geodesic Distance on Arbitrary-order Subdivision of Tree with Applications }}\\[12pt]
{\large Fei Ma$^{a,}$\footnote{~The author's E-mail: mafei123987@163.com. },\quad Xiaomin Wang$^{a,}$\footnote{~The author's E-mail: wmxwm0616@163.com.} \quad  and  \quad  Ping Wang$^{b,c,d,}$\footnote{~The author's E-mail: pwang@pku.edu.cn} }\\[6pt]
{\footnotesize $^{a}$ School of Electronics Engineering and Computer Science, Peking University, Beijing 100871, China\\
$^{b}$ National Engineering Research Center for Software Engineering, Peking University, Beijing, China\\
$^{c}$ School of Software and Microelectronics, Peking University, Beijing  102600, China\\
$^{d}$ Key Laboratory of High Confidence Software Technologies (PKU), Ministry of Education, Beijing, China}\\[12pt]
\end{center}

\begin{quote}
\textbf{Abstract:} The problem of how to estimate diffusion on a graph effectively is of importance both theoretically and practically. In this paper, we make use of two widely studied indices, geodesic distance and mean first-passage time ($MFPT$) for random walk, to consider such a problem on some treelike models of interest. To this end, we first introduce several types of operations, for instance, $m$th-order subdivision and ($1,m$)-star-fractal operation, to generate the potential candidate models. And then, we develop a class of novel techniques based on mapping for calculating the exact formulas of the both quantities above on our models. Compared to those previous tools including matrix-based methods for addressing the issue of this type, the techniques proposed here are more general mainly because we generalize the initial condition for creating these models. Meantime, in order to show applications of our techniques to many other treelike models, the two popularly discussed models, Cayley tree $C(t,n)$ and Exponential tree $\mathcal{T}(t;m)$, are chosen to serve as representatives. While their correspondence solutions have been reported, our techniques are more convenient to implement than those pre-existing ones and thus this certifies our statements strongly. Final, we distinguish the difference among results and attempt to make some heuristic explanations to the reasons why the phenomena appear.


\end{quote}

\vskip 1cm

\section{INTRODUCTION}

Complex network, also customarily called graph, has received extensively attention due to its own right in the recent years \cite{ Newman-2010}-\cite{Ulrich-2019}. A network in its simplest form is a collection of vertices connected together in pairs by edges. For instance, in a scientific citation network, a vertex represents one published paper and an edge is citation relation between a pair of vertices. Similarly, a vertex in the friend network is an individual and the friendship between a pair of individuals is denoted by an edge. It is these such networked models that help us to unveil some rules behind their topological structure including small-world property \cite{WS-1998} and scale-free feature \cite{BA-1999}. In addition, researchers are also interested in understanding the dynamics and functions appearing in these models. One of the best studied topics is to determine how fast information diffusion on complex networks is in part because of a great number of applications, such as, the study of transport-limited reactions \cite{O-B-2011}, target search \cite{Michael-2006}, disease spreading on relationship networks among individuals \cite{Jia-2018} as well navigation in spatial networks \cite{Wei-2014}, to name just a few.

In order to estimate the efficiency of information diffusion, there have been various measures proposed. The simplest of which is diameter $D$. For a given graph, its diameter is the maximum among geodesic distance of all pairs of vertices. Geodesic distance $d_{uv}$ of a pair of vertices $u$ and $v$, as well named shortest path distance, is the total number of edges of a shortest path whose endvertices are the both ones. In particular, it is using the diameter that almost all complex networks have been analytically proven to be small. Considering a graph $G(V,E)$ as a whole, one can also adopt the average geodesic distance $\langle \mathcal{S}\rangle$ defined here
\begin{equation}\label{Section-11}
\langle \mathcal{S}\rangle=\frac{\mathcal{S}}{|V|(|V|-1)/2}=\frac{\sum_{u,v\in V}d_{uv}}{|V|(|V|-1)/2}.
\end{equation}
There is no question that it is equivalent for a general graph to calculate the precise solutions for geodesic distance and its diameter, respectively. On the other hand, it is easy to obtain the diameters of some graphs with specific topological structure, in which case the closed-form formula of geodesic distance is too hard to solve. Therefore, it is an interesting challenge to find out some available manners in which one may deal with the problems of this type. Besides that, from the point of view of uncertainty, the random walk has proven with success in measuring the diffusion on graphs \cite{MaF-2019}-\cite{Tobias-2018}. The random-walk model considered here is a simple one. At each time step, the walker will jump at random to any of its neighbor set $N_{u}$ from its current position $u$ with equal probability, i.e.,

\begin{equation}\label{Section-12}
P_{u\rightarrow v}=\left\{\begin{aligned}&1/d_{_{u}} \qquad \text{if $u$ is adjacent to $v$},\\
&0 \qquad \qquad\text{otherwise.}
\end{aligned}\right.
\end{equation}
In other words,

$$fP(v)=\sum_{u,u\sim v}\frac{1}{d_{u}}f(u)$$
for any distribution $f:V\rightarrow \mathbf{R}$ with $\sum_{v}f(v)=1$ \cite{Chung-1998}. By analogy with geodesic distance $d_{uv}$, one usually focuses on the first-passage time for a walker starting from vertex $u$ to reach the destination $v$ for the first-passage time on graph, denoted by $F_{uv}$. As before, we can define a quantity referred to as the mean first-passage time ($MFPT$) in the following term
\begin{equation}\label{Section-13}
MFPT=\frac{\sum_{u\in V}\sum_{v\neq u,v\in V}F_{uv}}{|V|(|V|-1)}.
\end{equation}

In some extent, this certainly provides a helpful indicator for the efficiency of information diffusion on graphs. The goal of this paper is to primarily discuss the latter two measures on some trees of great interest as we will address later. As the most fundamental connected graph, tree has a great number of intriguing properties \cite{Bondy-2008}. Here, we make use of one of them to built up our main results. Such a property is that there is only a unique path for a given vertex pair.

The rest of this paper is organized as follows. We in Section 2 introduce some conventional notations in the jargon of graph theory \cite{Bondy-2008}. In addition, we also introduce some operations to construct some potential candidate models, for instance, $m$th-order subdivision and ($w,m$)-star-fractal operation. In principle, the operations proposed here can be applied to any graph and thus have been widely studied in a great range of science community in the past years. Our goal, however, is to focus on geodesic distance on the resulting trees based on the operations. To this end, in Section 3, we develop a series of novel techniques on the basis of mapping. By using the developed methods, we can calculate the exact solution for geodesic distance on treelike models in question. Compared to some previous tools for addressing such a problem, our techniques are more general and more convenient to carry out. To further verify the potential applications of our manners to many other treelike models, we generalize the construction of the both popularly studied trees, Cayley tree $C(t,n)$ and Exponential tree $\mathcal{T}(t;m)$, using two algorithms upon an arbitrary tree $\mathcal{T}$ in Section 4. In a special case, we obtain the precise formula of geodesic distance on Cayley tree $C(t,n)$ by virtue of its self-similarity. Nevertheless, such a method can not be adequately adopted when its seed may be an arbitrary tree. This suggests indirectly the universality of our techniques. The topic of Section 5 focuses mainly on studying random walks on the generated models. With the help of effective resistance of electrical networks, we capture the exact solution of mean first-passage time ($MFPT$) for a walker on each of the models under consideration. Armed with the consequences, we provide some discussions aiming at distinguishing similarities and differences among them in Section 6. We conclude that different operations have substantial influence on the $MFPT$ on the resulting graphs. Final, we close our work using a precise conclusion.

\section{DEFINITIONS AND NOTATIONS}
It is a convention in graph theory to denote a graph by $G(V,E)$ whose vertex set and edge set are $V$ and $E$, respectively, and hence vertex number (order) and edge number (size) are equal to $|V|$ and $|E|$ where symbol $|\;|$ indicates the cardinality of a set. At the same time, let notation $[1,n]$ be an integer set which consists certainly of those integers no greater than $n$ and no less than $1$. Refer to Ref.\cite{Bondy-2008} for more details.

\textbf{Definition 1} Given an arbitrary graph $G(V,E)$, if one inserts a new vertex $\omega$ to every edge $uv\in E$ then the resulting graph, denoted by $G'_{1}(V'_{1},E'_{1})$, is called a \textbf{\emph{first-order subdivision graph}} of the original graph. Put this another way, such a graph $G'_{1}(V'_{1},E'_{1})$ can be obtained from graph $G(V,E)$ by replacing every edge $uv\in E$ by a unique path $uwv$ with length two. Here we regard this operation implemented on edge as the \textbf{\emph{first-order subdivision }}vividly. Fig.1(b) shows such an operation on tree $\mathcal{T}$ on seven vertices plotted in Fig.1(a).

\begin{figure}
\centering
  \includegraphics[height=5cm]{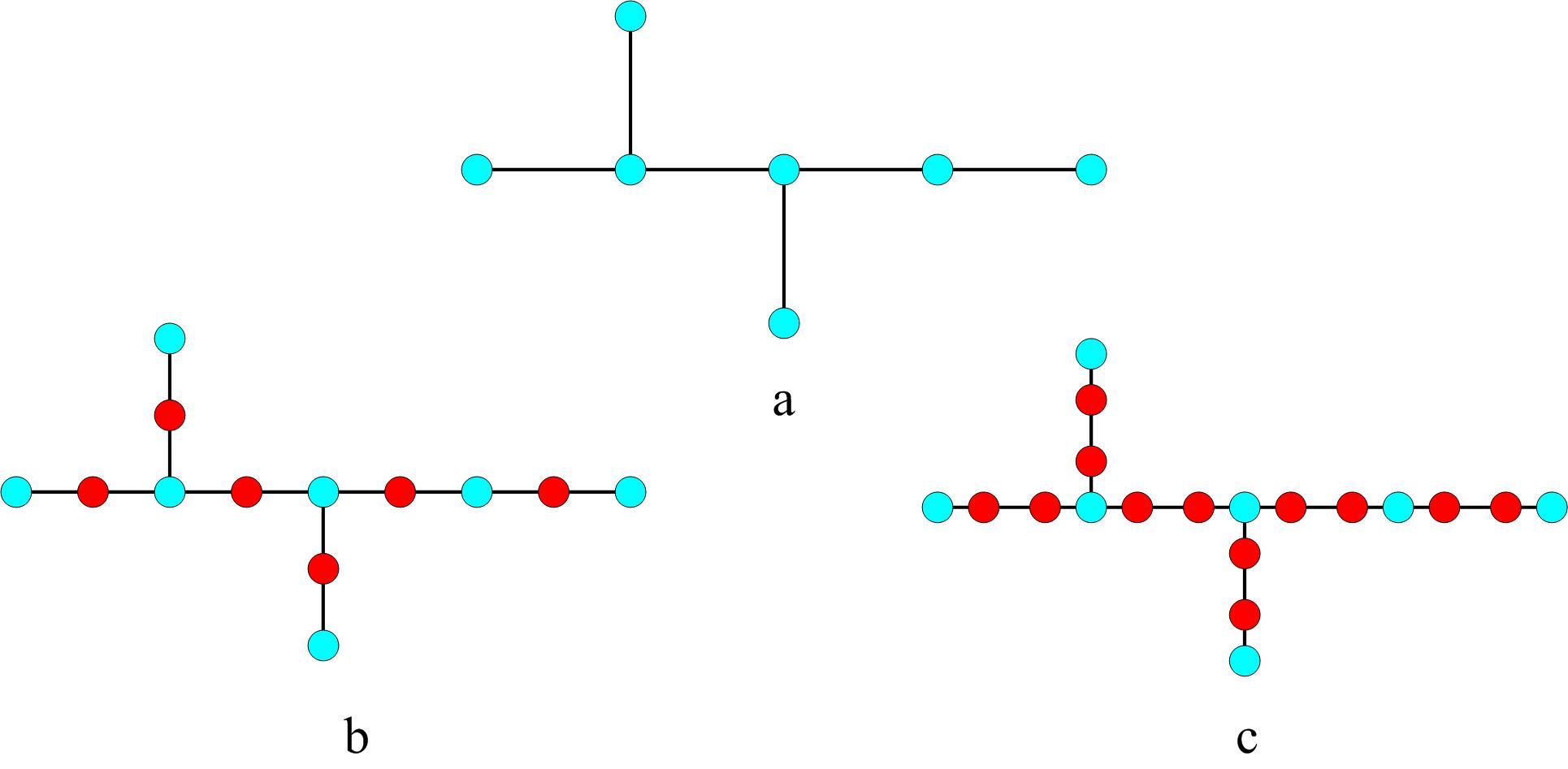}\\
{\small Fig.1. The diagrams of a tree $\mathcal{T}$ on seven vertices and its two types of subdivision trees. Panel (a) illustrates an original tree $\mathcal{T}$. The bottommost two show the first-order subdivision tree $\mathcal{T'}$ and $m$th-order subdivision tree $\mathcal{T}^{m}$ where $m=2$, respectively.}
\end{figure}

For a seminal graph $G(V,E)$, it can immediately see using Def.1 that its first-order subdivision graph $G'_{1}(V'_{1},E'_{1})$ holds on $|V'_{1}|=|V|+|E|$ and $|E'_{1}|=2|E|$. After $t$ time steps, by using an iterative method, we can obtain exactly the order $|V'_{t}|$ and size $|E'_{t}|$ of the first-order subdivision graph $G'_{t}(V'_{t},E'_{t})$ as follows
\begin{equation}\label{Section-2-D10}
|V'_{t}|=|V|+(2^{t}-1)|E|, \qquad |E'_{t}|=2^{t}|E|.
\end{equation}

Analogously, for a given graph $G(V,E)$, one would like to insert $m$ $(m\geq2)$ new vertices $\omega_{i}$ $(i\in [1,m])$ not a single vertex to every edge $uv\in E$, leading to the so-called \textbf{\emph{mth-order subdivision graph}}  $G^{m}_{1}(V^{m}_{1},E^{m}_{1})$ shown in Fig.1(c) where $m=2$. As before, such an operation is thought of as the \textbf{\emph{mth-order subdivision}} directly. After $t$ time steps, we are able to use the same method as developing Eq.(\ref{Section-2-D10}) to calculate the exact solutions for the order $|V^{m}_{t}|$ and size $|E^{m}_{t}|$ of the $m$th-order subdivision graph $G^{m}_{t}(V^{m}_{t},E^{m}_{t})$ in the following form
\begin{equation}\label{Section-2-D11}
|V^{m}_{t}|=|V|+[(m+1)^{t}-1]|E|, \qquad |E^{m}_{t}|=(m+1)^{t}|E|.
\end{equation}

\textbf{Definition 2} Given an arbitrary graph $G(V,E)$, if one not only inserts a new vertex $\omega$ to every edge $uv\in E$ but also connects $m$ other new vertices $\omega_{i}$ $(i\in [1,m])$ to this newly added vertex, then the resulting graph $G^{\star}_{1;1,m}(V^{\star}_{1;1,m},E^{\star}_{1;1,m})$ is a\textbf{ \emph{(1,m)-star-fractal graph} }of the original graph $G(V,E)$. As above, the ($1,m$)-star-fractal graph $G^{\star}_{1;1,m}(V^{\star}_{1;1,m},E^{\star}_{1;1,m})$ can be resulted from graph $G(V,E)$ by equivalently inserting a star with $m$ leaves to every $uv\in E$ and hence this operation is called\textbf{ \emph{(1,m)-star-fractal operation}}. As widely studied in a number of literature \cite{Peng-2018}, the well known \textbf{\emph{T-graph}} can be obviously induced as a special case of our ($1,m$)-star-fractal graph when the seminal graph $G(V,E)$ is defined as a single edge and parameter $m$ is supposed equal to $1$. An example as illustration of the ($1,m$)-star-fractal of tree $\mathcal{T}$ with seven vertices is shown in Fig.2(b) where $m=2$.

\begin{figure}
\centering
  \includegraphics[height=4cm]{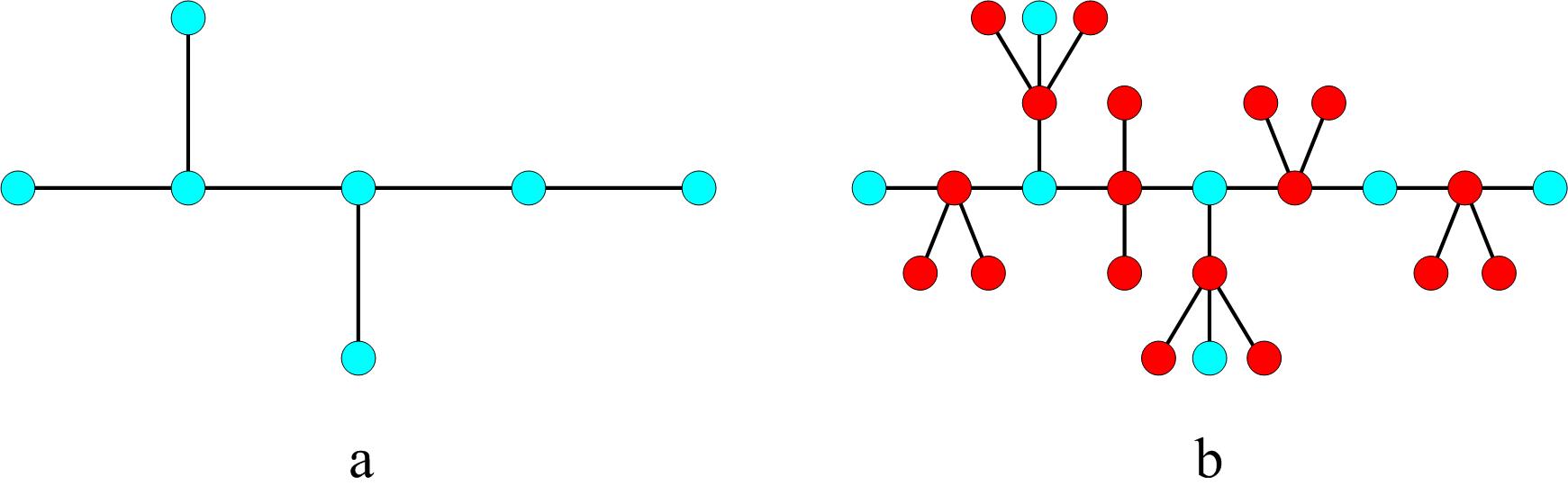}\\
{\small Fig.2. The diagrams of a tree $\mathcal{T}$ with seven vertices shown in panel (a) and its ($1,m$)-star-fractal tree $\mathcal{T}^{\star}_{1;1,m}$ where $m=2$ as plotted in panel (b).       }
\end{figure}

For a seminal graph $G(V,E)$, based on Def.2, one can easily understand that the ($1,m$)-star-fractal graph $G^{\star}_{1;1,m}(V^{\star}_{1;1,m},E^{\star}_{1;1,m})$ has $|V^{\star}_{1;1,m}|=|V|+(1+m)|E|$ vertices and $|E^{\star}_{1;1,m}|=(2+m)|E|$ edges. By the same as discussion in Def.1, for time step $t$, vertex number $|V^{\star}_{t;1,m}|$ and edge number $|E^{\star}_{t}|$ of the ($1,m$)-star-fractal graph $G^{\star}_{t;1,m}(V^{\star}_{t;1,m},E^{\star}_{t;1,m})$, respectively, follow
\begin{equation}\label{Section-2-D20}
|V^{\star}_{t;1,m}|=|V|+[(2+m)^{t}-1]|E|, \qquad |E^{\star}_{t;1,m}|=(2+m)^{t}|E|.
\end{equation}

Taking into account the $w$th-order subdivision in Def.1, for an arbitrary graph $G(V,E)$, one may attempt to connect $m$ other new vertices $\omega_{i}$ $(i \in [1,m])$ to each of these newly added vertices $u_{i}$ $(i\in [1,w])$ of the $w$th-order subdivision graph $G^{w}_{1}(V^{w}_{1},E^{w}_{1})$, resulting in the so-called\textbf{ \emph{(w,m)-star-fractal graph} }$G^{\star}_{1;w,m}(V^{\star}_{1;w,m},E^{\star}_{1;w,m})$. In fact, such an operation can be equally achieved by directly adding $w$ stars with $m$ leaves to each edge $uv$ of graph $G(V,E)$ and thus is viewed as the\textbf{ \emph{(w,m)-star-fractal operation}}, seeing Fig.3. As mentioned in Def.1, after $t$ time steps, we make use of an iterative method to calculate the precise solutions for the order $|V^{\star}_{t;w,m}|$ and size $|E^{\star}_{t;w,m}|$ of the ($w,m$)-star-fractal graph $G^{\star}_{t;w,m}(V^{\star}_{t;w,m},E^{\star}_{t;w,m})$ as below
\begin{equation}\label{Section-2-D22}
|V^{\star}_{t;w,m}|=|V|+([(w+1)m+1]^{t}-1)|E|, \qquad |E^{\star}_{t;w,m}|=[(w+1)m+1]^{t}|E|.
\end{equation}

With the definition of average degree $\langle k\rangle=2|E|/|V|$, we remark that for an arbitrary graph $G(V,E)$, after $t$ time steps, the four resulting graphs above all have an identical average degree equal to $2$ in the limit of large graph size. Nonetheless, there exist significant differences among these models as we will show later.

\begin{figure}
\centering
  \includegraphics[height=4.5cm]{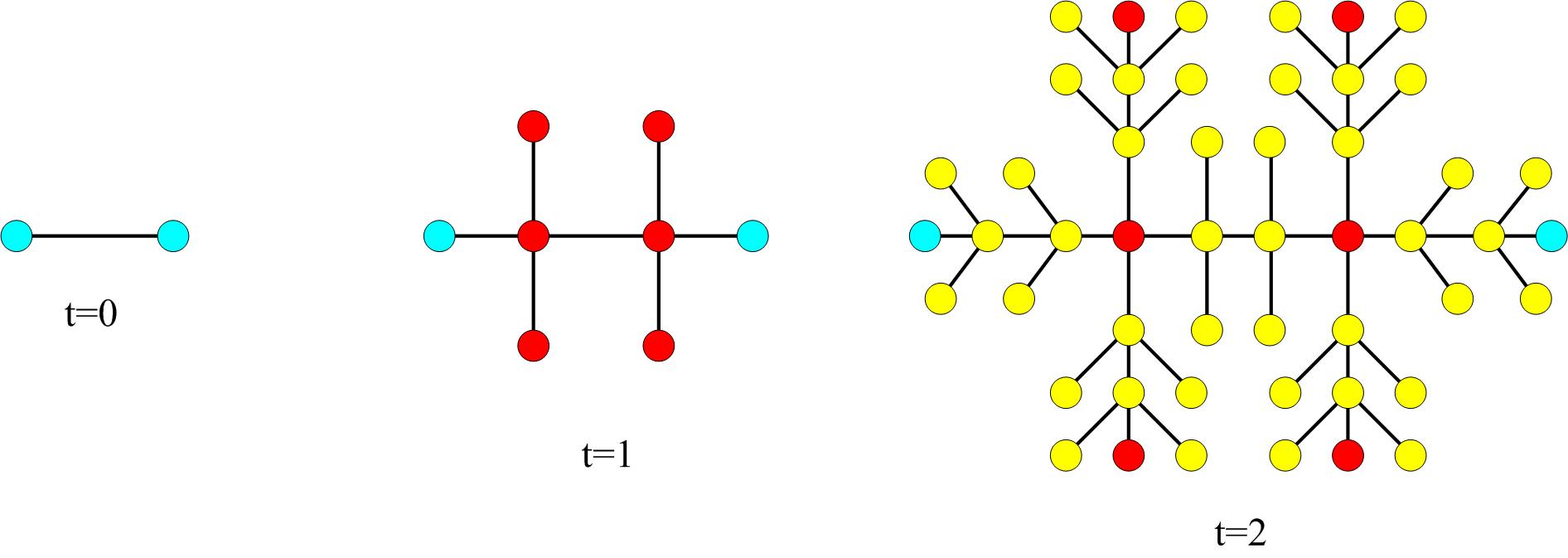}\\
{\small Fig.3. The diagram of the ($w,m$)-star-fractal operation on an edge where $w=m=2$.       }
\end{figure}

\textbf{Definition 3} Given two nonempty sets $X$ and $Y$, one can take a \textbf{\emph{mapping}} $f$ between $X$ and $Y$ such that for a provided element $x$ of set $X$ there must be a unique element $y$ in set $Y$ satisfying $f(x)=y$. In the language of mathematics, this can be simply written as

$$\forall x\in X, \quad \exists y\in Y, \quad s.t.,\quad f:\; x\mapsto \; y$$
where the image set of set $X$ may be considered as $f_{X\mapsto Y}=\{y|f(x)=y, \; y\in Y\}$. Here we are particularly interested in the following two types of mappings between sets $X$ and $Y$.

\emph{case 1} If both $|X|>|Y|$ and $|f_{X\mapsto Y}|=|Y|$ are true, then this mapping $f$ is considered \textbf{\emph{surjection}}. Besides, for each element $y$ of image set $Y$, there indeed exist $n$ distinct pre-images $x_{i}$ ($i\in [1,n]$) in set $X$, i.e., $f^{-1}(y)=x_{i}$, then the surjection $f$ is considered \textbf{\emph{n-regular}}. It is clear to the eye that $|X|=n |Y|$ is true when the surjection $f$ in question is $n$-regular.

\emph{case 2} If the surjection $f$ under consideration holds $|X|=|Y|$ then it can be thought of as a \textbf{\emph{bijection}}, also called one-one mapping.

Meanwhile, the compound mapping between mappings $f$ and $g$ can be expressed as $f\circ g$ in the mathematical form.

As stated above, for a given graph $G(V,E)$, Eqs.(\ref{Section-2-D10})-(\ref{Section-2-D22}) hold truth. The goal of this paper, however, focuses mainly on many discussions about geodesic distance on treelike models where the seminal graph is an arbitrary tree $\mathcal{T}$. In addition, we will take vertex pair $<u,v>$ to represent the geodesic distance of the same pair of vertices    because tree is the simplest connected graph. In other words, there is only a unique path in a tree which connects vertex $u$ to vertex $v$. Now, let us divert our insight to problems.

\section{MAIN RESULTS AND APPLICATIONS}

\textbf{Lemma 1} For a path $\mathcal{P}$ with $A$ vertices ($A\geq2$), its corresponding geodesic distance $\mathcal{S}_{p}$ obeys

\begin{equation}\label{Section-3-l1}
\mathcal{S}_{p}=\sum_{i=1}^{A-1}\sum_{j=1}^{A-i}j.
\end{equation}

This is an obvious consequence in an enumerative method and hence we omit its proof.

\textbf{Lemma 2} For a tree $\mathcal{T}$ on $n$ vertices, the exact solution for geodesic distance on its first-order subdivision tree $\mathcal{T'}(1)$ is
\begin{equation}\label{Section-3-l2}
\mathcal{S}_{\mathcal{T'}(1)}=8\mathcal{S}-2n(n-1).
\end{equation}
where $\mathcal{S}$ is geodesic distance on tree $\mathcal{T}$.

\textbf{Proof} Suppose that there exists an arbitrary tree $\mathcal{T}$ with $n$ vertices and its geodesic distance is measured equal to $\mathcal{S}$. After applying the first-order subdivision to each edge of tree $\mathcal{T}$, the first-order subdivision tree $\mathcal{T'}(1)$ will consist of two different classes of vertices, without loss of generality, which are grouped into two disjoint vertex sets $X'$ composed of the old $n$ vertices and $Y'$ containing the young $n-1$ ones. In order to completely calculate geodesic distance $\mathcal{S'}$, it is necessary to compute three types of geodesic distances, one for vertex pairs $<x'_{i},x'_{j}>$ of set $X'$, one for vertex pairs $<y'_{i},y'_{j}>$ of set $Y'$ along with the latter for vertex pairs $<x'_{i},y'_{j}>$ between set $X'$ and set $Y'$. To this end, we will in turn accomplish these calculations.

\emph{Case 1.1} For a given vertex pair $<x'_{i},x'_{j}>$ of set $X'$, its geodesic distance is increased by a factor of $2$ according to the function of the first-order subdivision. Hence, one can have
\begin{equation}\label{Section-3-l2-1}
\mathcal{S'}_{1}(1)=2\mathcal{S}
\end{equation}
where $\mathcal{S'}_{1}(1)$ is the sum of geodesic distances over all possible vertex pairs belonging to set $X'$. The remaining tasks is to look for a reliable relation connecting the next equations to Eq.(\ref{Section-3-l2-1}) by establishing reasonable mappings between many other types of vertex pairs and the ones mentioned here.

\emph{Case 1.2} After applying the first-order subdivision to each edge, it is natural for a given path $<x'_{i},x'_{j}>$ whose two vertices are chosen from set $X'$ to insert $|j-i|$ new vertices labeled successively $y'_{i},y'_{i+1},...,y'_{j-1}$ (or equivalently marked $y'_{i+1},...,y'_{j-1},y'_{j}$). By using this consequence, we may introduce a bijection $f_{1}$ between set $Y'$ and set $X'$ such that $f_{1}(<y'_{i},y'_{j}>)=<x'_{i},x'_{j+1}>$ is true. Keep this in mind, one can expand a path $\mathcal{P}_{y'_{i}y'_{j}}$ as a unique path $\mathcal{P}_{x'_{i}x'_{j+1}}$ by connecting vertices $y'_{i}$ and $y'_{j}$ to vertices $x'_{i}$ and $x'_{j+1}$, respectively. Thus, we write
\begin{equation}\label{Section-3-l2-2}
\mathcal{S'}_{1}(2)=\mathcal{S'}_{1}(1)-2\left(\frac{n(n-1)}{2}\right)
\end{equation}
in which $\mathcal{S'}_{1}(2)$ is geodesic distance over all possible vertex pairs in set $Y'$.

\emph{Case 1.3} The task to answer is to capture the formula of geodesic distance for all possible vertex pairs $<x'_{i},y'_{j}>$ whose vertices are chosen from different sets, set $X'$ and set $Y'$. Along the research line of \emph{case 1.2}, there must be a unique path $\mathcal{P}_{x'_{i}y'_{j}}$ which can be reduced to another path $\mathcal{P}_{y'_{i}y'_{j}}$ by removing an additional edge $x'_{i}y'_{i}$ under a similar mapping $f_{2}$ to mapping $f_{1}$ in \emph{case 1.2}. Therefore, we may take a mapping $f_{2}$ which maps from vertex pair $<x'_{i},y'_{j}>$ to $<y'_{i},y'_{j}>$. Such a mapping is not bijection rather than surjection. To see why this is so, let we select a path $\mathcal{P}_{y'_{i}x'_{j+1}}$ and then have a choice to obtain path $\mathcal{P}_{y'_{i}y'_{j}}$ by deleting an edge $y'_{j}x'_{j+1}$. This suggests the generated mapping  $f_{2}$ is $2$-regular. Equivalently, there should be two distinct pre-images under surjection $f_{2}$, i.e., $<x'_{i},y'_{j}>$ and $<y'_{i},x'_{j+1}>$, such that $f^{-1}_{2}(<y'_{i},y'_{j}>)=<x'_{i},y'_{j}>=<y'_{i},x'_{j}>$. Armed with the both cases, we can obtain the exact solution for geodesic distance $\mathcal{S'}_{1}(3)$ over all possible vertex pairs $<x'_{i},y'_{j}>$ and $<y'_{i},x'_{j+1}>$ as follows
\begin{equation}\label{Section-3-l2-3}
\mathcal{S'}_{1}(3)=2\mathcal{S'}_{1}(2)+n(n-1).
\end{equation}

With the results in \emph{cases 1.1-1.3}, we have exhaustively enumerated all possible vertex pairs. Therefore, Eqs.(\ref{Section-3-l2-1})-(\ref{Section-3-l2-3}) together yields an exact formula for $\mathcal{S'}_{1}$ in complete agreement with that of Eq.(\ref{Section-3-l2}). This completes our proof.

By using the consequence of Eq.(\ref{Section-3-l2}), we can capture the solution for geodesic distance $\mathcal{S'}_{t}$ on the first-order subdivision tree $\mathcal{T'}(t)$ by adopting some simple arithmetics.

\emph{Corollary 1} For a tree $\mathcal{T}$ on $n$ vertices, the exact solution for geodesic distance on its first-order subdivision tree $\mathcal{T'}(t)$ is
\begin{equation}\label{Section-3-c1}
\mathcal{S}_{\mathcal{T'}(t)}=8^{t}\mathcal{S}-\frac{1}{3}(2^{3t}-2^{t})(n-1)+(2^{2t-1}-2^{3t-1})(n-1)^{2}.
\end{equation}

If one is interested in selecting a single edge $\mathcal{E}_{uv}$ as a seed, then the first-order subdivision ``edge'' $\mathcal{E'}(t)$ is a path with $2^{t+1}+1$ vertices. Combining the results from lemma 1 and lemma 2, we can arrive at the following corollary 2.

\emph{Corollary 2} For an edge $\mathcal{E}$, the exact solution for geodesic distance on its first-order subdivision ``edge" $\mathcal{E'}(t)$ is
\begin{equation}\label{Section-3-c2}
\mathcal{S}_{\mathcal{E'}(t)}=\sum_{i=1}^{2^{t}}\sum_{j=1}^{2^{t}+1-i}j=\frac{(2^{t}+1)2^{t}(2^{t-1}+1)}{3}.
\end{equation}

Here, Eq.(\ref{Section-3-c2}) provides a proof of a combinatorial identity in a novel manner where parameter $A$ can be expressed as $A=\sum_{i=0}^{t}2^{i}+1$. Nevertheless, Eq.(\ref{Section-3-l1}) indicates a general case in which parameter $A$ may be equivalent to any integer greater than $2$. To generalize our result addressed in Eq.(\ref{Section-3-c1}) to accomplish our desired consequences, we need to introduce a fact from number theory, i.e.,
$$A=\sum_{i=0}^{n}a^{i}-1,\qquad \forall A, \forall a\in \mathrm{N}_{+}.$$

Since then, we can choose an available basis $a$ to expand parameter $A$ and further calculate the exact formula of  Eq.(\ref{Section-3-l1}). This can be achieved using the results in the below theorem as we will show shortly.

\textbf{Theorem 1} For a tree $\mathcal{T}$ on $n$ vertices, the exact solution for geodesic distance on its $m$th-order subdivision tree $\mathcal{T}^{m}(1)$ is
\begin{equation}\label{Section-3-t1}
\mathcal{S}_{\mathcal{T}^{m}(1)}=(m+1)^{3}\mathcal{S}-\frac{2n^{2}\sum_{i=1}^{m}i^{3}}{m}+2n\sum_{i=1}^{m}i^{2}-\frac{m+1}{2m-1}\sum_{i=1}^{m-1}i^{2}.
\end{equation}

\textbf{Proof} Consider a tree $\mathcal{T}$ on $n$ vertices, let us suppose its geodesic distance to be $\mathcal{S}$. As before, the $m$th-order subdivision tree $\mathcal{T}^{m}(1)$ also consists certainly of two sets, denoted by $X'$ and $Y'$. Vertices of tree $\mathcal{T}$ are grouped into set $X'$ and set $Y'$ is composed of all newly inserted vertices whose cardinality is $(n-1)m$. With the similar marks, for an edge $x'_{i}x'_{i+1}$, we label consecutively the $m$ young vertices inserted to this edge $y'_{i;1}, y'_{i;2},...,y'_{i;m}$. With the lights shed by the development of lemma 2, we still make use of classification method to determine contributions from different vertex pairs to geodesic distance $\mathcal{S}_{\mathcal{T}^{m}(1)}$ of tree $\mathcal{T}^{m}(1)$.

\emph{Case 2.1} For a pair of vertices $x'_{i}$ and $x'_{j}$, its geodesic distance is increased by a factor of $m+1$. Hence, the sum of geodesic distances over all possible vertex pairs of this type obeys

\begin{equation}\label{Section-3-t1-1}
\mathcal{S}^{(m)}(1)=(m+1)\mathcal{S}.
\end{equation}

\emph{Case 2.2} As known, there are $m$ new vertices embedded to each edge $uv$ of tree $\mathcal{T}$ by the $m$th-order subdivision. If we view each edge $uv$ of tree $\mathcal{T}$ as an ingredient, then for an arbitrary vertex pair $<y'_{u;i}, y'_{u;j}>$ $(i,j\in[1,m])$ its geodesic distance is equal to $|j-i|$. In addition, we should discuss geodesic distance over vertex pair $<y'_{s;i}, y'_{t;j}>$ $(i,j\in[1,m])$ whose endvertices may be chosen from two different ingredients, in which case we need to take similar measure to that in \emph{case 1.2} to build up a mapping $g_{1}$. Such a candidate $g_{1}$ can map vertex pair $<y'_{s;i}, y'_{t;j}>$ to vertex pair $<x'_{s}, x'_{t}>$, i.e., $g_{1}(<y'_{s;i}, y'_{t;j}>)=<x'_{s}, x'_{t}>$. In fact, this mapping turns out to be $m^{2}$-regular. Based on the two cases above, we can say
\begin{equation}\label{Section-3-t1-2}
\mathcal{S}^{(m)}(2)=m^{2}\mathcal{S}^{(m)}(1)+(n-1)\sum_{i=1}^{m-1}\sum_{j=1}^{m-i}j-\frac{n(n-1)(m+1)m^{2}}{2}.
\end{equation}
where $\mathcal{S}^{(m)}(2)$ represents the sum of geodesic distance over all possible vertex pairs in which two endvertices come from young vertex set $Y'$.

\emph{Case 2.3} Now we will measure the geodesic distance over any pair of vertices $y'_{s;i}$ $(i\in[1,m])$ and $x'_{t}$ $(s,t\in[1,n])$. Analogously, a mapping $g_{2}$ should be generated to connect some known results to current problems to answer. Without loss of generality, we assume that the vertex $y'_{s;i}$ is one of new vertices inserted to an ingredient $x'_{s}x'_{s+1}$. Under mapping $g_{2}$, we require that the path $\mathcal{P}_{y'_{s;i}x'_{t}}$ be reduced as a path $\mathcal{P}_{x'_{s+1}x'_{t}}$ whose geodesic distance has be obtained from Eq.(\ref{Section-3-t1-1}). Meantime, it is not hard to prove such a mapping $g_{2}$ to be $2m$-regular. One reason for this is that the vertex pairs $<x'_{s+1},y'_{t-1;i}>$ are also pre-images of vertex pair $<x'_{s+1},x'_{t}>$ under mapping $g_{2}$. Through these statements above, the geodesic distance $\mathcal{S}^{(m)}(3)$ over all possible vertex pairs of this type should follow

\begin{equation}\label{Section-3-t1-3}
\mathcal{S}^{(m)}(3)=2n\mathcal{S}^{(m)}(1)-n(n-1)\sum_{i=1}^{m}i.
\end{equation}

By using some simple and necessary arithmetics, Eqs.(\ref{Section-3-t1-1})-(\ref{Section-3-t1-3}) can yield a desirable consequence which is consistent with that of Eq.(\ref{Section-3-t1}). We complete the proof of theorem 1.

After $t$ time steps, we can immediately arrive at the following corollary.

\emph{Corollary 3} For a tree $\mathcal{T}$ on $n$ vertices, the exact solution for geodesic distance on its $m$th-order subdivision tree $\mathcal{T}^{m}(t)$ is
\begin{equation}\label{Section-3-c3}
\begin{aligned}\mathcal{S}_{\mathcal{T}^{m}(t)}&=(m+1)^{3t}\mathcal{S}-\frac{2(n-1)^{2}\sum_{i=1}^{m}i^{3}\sum_{j=0}^{t-1}(m+1)^{2(t-1)+j}}{m}-\frac{4(n-1)\sum_{i=1}^{m}i^{3}\sum_{j=0}^{t-1}(m+1)^{2(t+j)-1}}{m}\\
&-\frac{2\sum_{i=1}^{m}i^{3}\sum_{j=0}^{t-1}(m+1)^{3j}}{m}+2(n-1)\sum_{i=1}^{m}i^{2}\sum_{j=0}^{t-1}(m+1)^{t+2j-1}-\frac{m+1}{2n-1}\sum_{i=1}^{m-1}i^{2}\sum_{j=0}^{t-1}(m+1)^{3j}\\
&+2\sum_{i=1}^{m}i^{2}\sum_{j=0}^{t-1}(m+1)^{3j}.
\end{aligned}
\end{equation}

It is worthy noting that we indeed finish a more general combinatorial identity in terms of Eq.(\ref{Section-3-l1}) and Eq.(\ref{Section-3-c3}) in comparison with Eq.(\ref{Section-3-c2}). Nonetheless, if the original tree degenerates an edge $\mathcal{E}$, the geodesic distance $\mathcal{S}_{\mathcal{E}^{m}(t)}$ of the $m$th-order subdivision ``edge" $\mathcal{E}^{m}(t)$ can be obtained by using either Eq.(\ref{Section-3-l1}) where $A=(m+1)^{t}+1$ or Eq.(\ref{Section-3-c3}) where $\mathcal{S}=1$ and $n=2$. From the complexity point of view, Eq.(\ref{Section-3-l1}) should be adequately employed in such a situation.

By far, we have completely calculated the formulas of geodesic distance on $m$th-order subdivision tree $\mathcal{T}^{m}(t)$ whose seed is an arbitrary tree $\mathcal{T}$ we please. From now on, we turn to discuss a family of treelike models of interest which are produced on the basis of the ($w,m$)-star-fractal operation. As a start point, we directly bring a lemma and a corollary with omitting calculations. Interested reader is referred to see our previous work \cite{X-M-2019}.

\textbf{Lemma 3} For a tree $\mathcal{T}$ on $n$ vertices, the exact solution for geodesic distance on its ($1,m$)-star-fractal tree $\mathcal{T}^{\star}(1;1,m)$ is

\begin{equation}\label{Section-3-l3}
\mathcal{S}_{\mathcal{T}^{\star}(1;1,m)}=2(m+2)^{2}\mathcal{S}-(m+2)(n-1)(m+n).
\end{equation}

\emph{Corollary 4} For a tree $\mathcal{T}$ on $n$ vertices, the exact solution for geodesic distance on its ($1,m$)-star-fractal tree $\mathcal{T}^{\star}(t;1,m)$ is

\begin{equation}\label{Section-3-c4}
\begin{aligned}\mathcal{S}_{\mathcal{T}^{\star}(t;1,m)}&=2^{t}(m+2)^{2t}\mathcal{S}-(2^{t}-1)(m+2)^{2t-1}(n^{2}-2n-1)\\
&-\frac{(m+1)(n-1)}{2}\times\frac{2^{t+1}(m+2)^{2t}-2(m+2)^{t}}{2m+3}.
\end{aligned}
\end{equation}

Now, let us go forth into the discussion about geodesic distance on the ($w,m$)-star-fractal tree $\mathcal{T}^{\star}(1;w,m)$.

\textbf{Theorem 2} For a tree $\mathcal{T}$ on $n$ vertices, the exact solution for geodesic distance on its ($w,m$)-star-fractal tree $\mathcal{T}^{\star}(1;w,m)$ is

\begin{equation}\label{Section-3-t2}
\begin{aligned}\mathcal{S}_{\mathcal{T}^{\star}(1;w,m)}&=(w+1)[w(m+1)+1]^{2}\mathcal{S}-\frac{wn^{2}}{2}[(m+1)^{2}w^{2}-(m-2)(m+1)w-(m-1)]\\
&+\frac{wn}{6}[4(m+1)^{2}w^{2}-3(m+1)(3m-2)w-(m^{2}+11m+2)]\\
&-\frac{w}{6}[(m+1)^{2}w^{2}-6m(m+1)w-(m^{2}+8m+1)]\\
&=\Psi_{1}(m,w)\mathcal{S}-\Psi_{2}(m,w)n^{2}+\Psi_{3}(m,w)n-\Psi_{4}(m,w).
\end{aligned}
\end{equation}

\textbf{Proof} As the ($w,m$)-star-fractal operation shows, there are $w$ stars with $m$ leaf vertices each added to each edge of the seminal tree $\mathcal{T}$. Thus, we group vertices of ($w,m$)-star-fractal tree $\mathcal{T}^{\star}(1;w,m)$ into two sets, set $X'$ containing all vertices of tree $\mathcal{T}$ and the other $Y'$ composed entirely of those newly generated vertices. It turns out that the machinery developed above can be used to help us to complete the proof of theorem 2. The only thing that needs to be changed is to divide the vertex set $Y'$ into many smaller sets. Equivalently, the set $Y'$ is a union of its several disjoints subsets. Here we denote by the subset $Y'_{1}$ a collection whose elements are selected from the central vertex of all stars. The other subset $Y'_{2}$ is made of leaf vertices of each star. In principle, our calculations have three parts. The first part represents the sum of geodesic distances over arbitrary pair of vertices in set $X'$, the second is equal to the sum of geodesic distances over any pair of vertices in set $Y'$, and the reminder is the sum of geodesic distance over all possible crossing vertex pairs. As we will shown shortly, the complete proof is in fact constituted of seven stages. Now, let us introduce them in stages.

\emph{Case 3.1} In terms of the characteristic of the ($w,m$)-star-fractal operation, the addition of the central vertex of each star will have influence on the geodesic distance over vertex pair of set $X'$. Since $w$ stars are inserted to each edge of tree $\mathcal{T}$, the geodesic distance sum $\mathcal{S}^{(w)}_{m}(1)$ over all old vertex pairs satisfies

\begin{equation}\label{Section-3-t2-1}
\mathcal{S}^{(w)}_{m}(1)=(w+1)\mathcal{S}.
\end{equation}

\emph{Case 3.2} Between the $m$th-order subdivision and the ($w,m$)-star-fractal operation, the former can be easily reduced as a special case of the latter by deleting all leaf vertices of each star. This means that for arbitrary vertex pairs of subset $Y'_{1}$ the geodesic distance $\mathcal{S}^{(w)}_{m}(2)$ follows Eq.(\ref{Section-3-l2-2}) after replacing parameter $m$ by $w$.

\begin{equation}\label{Section-3-t2-2}
\mathcal{S}^{(w)}_{m}(2)=w^{2}\mathcal{S}^{(w)}_{m}(1)+(n-1)\sum_{i=1}^{w-1}\sum_{j=1}^{w-i}j-\frac{n(n-1)(w+1)w^{2}}{2}.
\end{equation}

\emph{Case 3.3} For each of newly created stars, there are $m$ leaf vertices attached to its central vertex and then the sum of geodesic distances over any pair leaf vertices equals $2[m(m-1)/2]$. Therefore, the geodesic distance sum $\mathcal{S}^{(w)}_{m}(3)$ over stars of such type is

\begin{equation}\label{Section-3-t2-3}
\mathcal{S}^{(w)}_{m}(3)=wm(m-1)(n-1).
\end{equation}

\emph{Case 3.4} By analogy with the trick in case 2.3, we take advantage of a mapping $h_{1}$ to connect the sum $\mathcal{S}^{(w)}_{m}(4)$ of geodesic distance over the first type of crossing vertex pairs to that of Eq.(\ref{Section-3-t2-1}). Here the first type of crossing vertex pairs indicates one of which is chosen at random from vertex set $Y'_{1}$ and the other chosen randomly from vertex set $X'$. In practice, such a formula has been obtained as shown in Eq.(\ref{Section-3-t1-3}). The only requirement is the exchange between current parameters and the previous and we omit the explanation in more detail.

\begin{equation}\label{Section-3-t2-4}
\mathcal{S}^{(w)}_{m}(4)=2n\mathcal{S}^{(w)}_{m}(1)-n(n-1)\sum_{i=1}^{w}i.
\end{equation}

\emph{Case 3.5} Following the discussion of in case 3.4, we now compute the geodesic distance $\mathcal{S}^{(w)}_{m}(5)$ of the second type of crossing vertex pairs. Given such a pair of vertices, one vertex is selected from set $Y'_{2}$ and the other is still in set $X'$. Obviously, vertex in set $Y'_{2}$ is connected by an edge to the central vertex of star to which it belongs. Combining the fact with the known result in Eq.(\ref{Section-3-t2-4}), we may create a mapping $h_{2}$ which maps the vertex pair in question to a unique path whose geodesic distance is discussed in case 3.4 by removing that edge linking leaf vertex to the central vertex of star to which it belongs. With the help of compound mapping $h_{2}\circ h_{1}$, we can have

\begin{equation}\label{Section-3-t2-5}
\mathcal{S}^{(w)}_{m}(5)=m\mathcal{S}^{(w)}_{m}(4)+wmn(n-1).
\end{equation}

\emph{Case 3.6} Now let us pay attention on the computation of geodesic distance of an arbitrary pair of vertices. On the one hand, the both vertices of the vertex pair of this type are in set $Y'_{2}$. On the other hand, the both do belong to the same star. Without loss of generality, if we denote by $y'_{i,j}$ a vertex that is a leaf labeled $j$ attached to the central vertex $y'_{i}$ of a star labeled $i$ added to an edge $uv$ of the seminal tree $\mathcal{T}$, then the notation $y'_{i,j}y'_{\alpha,\beta}$ can be regarded as a couple of vertices under consideration. To accomplish the sum of geodesic distance over all possible pair of vertices in question, we need to adopt the thought reflected by the proceeding cases, namely, the introduction of new mapping $h_{3}$. The performance of mapping $h_{3}$ is in spirit similar to the previous mapping $h_{1}$. For a pair of vertices $y'_{i,j}y'_{\alpha,\beta}$ in which the central vertices $y'_{i}$ and $y'_{\alpha}$ are inserted to edges $uv$ and $st$, respectively, this $h_{3}$ maps path $<y'_{i,j},y'_{\alpha,\beta}>$ into path $<y'_{i},y'_{\alpha}>$ by the removal of two edges $y'_{j}y'_{i}$ and $y'_{\beta}y'_{\alpha}$. The geodesic distance of the latter has been shown in Eq.(\ref{Section-3-t2-2}). Considering that the central vertex of each star is attached to $m$ distinct leaf vertices, it is not hard to see that mapping $h_{3}$ is in fact $m^{2}-$regular. Based on the description stated here and Eq.(\ref{Section-3-t2-2}), the precise solution for the geodesic distance $\mathcal{S}^{(w)}_{m}(6)$ of vertex pairs $<y'_{i,j},y'_{\alpha,\beta}>$ of such type obeys

\begin{equation}\label{Section-3-t2-6}
\begin{aligned}\mathcal{S}^{(w)}_{m}(6)&=m^{2}\mathcal{S}^{(w)}_{m}(2)+\sum_{i=1}^{w(n-1)-1}2m^{2}[w(n-1)-i]\\
&=(w+1)(wm)^{2}\mathcal{S}-\frac{(w+1)(wm)^{2}}{2}n^{2}+\frac{w(w+1)m^{2}}{6}[n(4w-1)-(w-1)]\\
&+(wm)^{2}(n^{2}-3n+2)+w(w-1)m^{2}(n-1).
\end{aligned}
\end{equation}

\emph{Case 3.7} So far, we have finished the calculation of geodesic distance over vertex pairs whose two endvertices are either  in set $X'$, or in $Y'_{i}$, even one in $X'$ and the other in $Y'$. The rest of our task is to determine the formula for the geodesic distance over any pair of vertices $<y'_{i},y'_{\alpha,\beta}>$ where vertex $y'_{i}$ is in set $Y'_{1}$ and vertex $y'_{\alpha,\beta}$ is a leaf belonging to set $Y'_{2}$. As before, we still want to employ the mapping-based technique to capture the such calculation. To be specific, given a vertex pair $<y'_{i},y'_{\alpha,\beta}>$, one may build up a $2m-$regular mapping $h_{4}$ such that $h_{4}(<y'_{i},y'_{\alpha,\beta}>)=<y'_{i},y'_{\alpha}>$. It is easy to achieve such a mapping because the path $\mathcal{P}_{y'_{i},y'_{\alpha,\beta}}$ can be reduced as a distinct path $\mathcal{P}_{y'_{i},y'_{\alpha}}$ by the deletion of an edge $y'_{\alpha}y'_{\alpha,\beta}$. The $2m-$regularity of mapping $h_{4}$ is in view of the central vertex of each star attached to $m$ leaf vertices. Once mapping $h_{4}$ is fixed, one can switch the computation of the sum $\mathcal{S}^{(w)}_{m}(7)$ of geodesic distance of vertex pairs in question into the known result in Eq.(\ref{Section-3-t2-2}) by some simple arithmetic, as follows

\begin{equation}\label{Section-3-t2-7}
\begin{aligned}
\mathcal{S}^{(w)}_{m}(7)&=2m\mathcal{S}^{(w)}_{m}(3)+m[w(n-1)]^{2}\\
&=2m(w+1)w^{2}\mathcal{S}-m(w+1)w^{2}n^{2}+\frac{w(w+1)(4w-1)m}{3}n\\
&-\frac{(w-1)w(w+1)m}{3}+mw^{2}n^{2}-2mw^{2}n+mw^{2}.
\end{aligned}
\end{equation}

Till now, we have exhaustively enumerated the geodesic distance over all possible vertex pairs on the ($w,m$)-star-fractal tree $\mathcal{T}^{\star}(1;w,m)$. Armed with the results above, we are able to understand

\begin{equation}\label{Section-3-t2-8}
\begin{aligned}\mathcal{S}^{(w)}_{m}&=\sum_{i=1}^{7}\mathcal{S}^{(w)}_{m}(i)=[(w+1)^{3}+wm(w+1)(wm+2w+2)]\mathcal{S}\\
&-\left[\frac{w(w+1)^{2}}{2}+\frac{wm(w-1)(wm+1)}{2}+mw^{3}\right]n^{2}+\left[\frac{w(w+1)(2w+1)}{3}+\frac{w(4w^{2}-9w-7)m^{2}}{6}\right]n\\
&+\left[wm(m-1)+\frac{wm(w-1)(8w+5)}{2}\right]n-\left[\frac{(w-7)w(w+1)m^{2}}{6}+\frac{wm(w^{2}-3w-1)}{3}\right]\\
&-\left[wm(m-1)+\frac{(w-1)w(w+1)}{6}\right].
\end{aligned}
\end{equation}
as desired. By using some simple arithmetics, Eq.(\ref{Section-3-t2-8}) will be completely consistent with Eq.(\ref{Section-3-t2}) and thus this completes the proof of theorem 2.

As an immediate consequence of theorem 2, if one would like to apply the ($w,m$)-star-fractal operation to a given seed tree $\mathcal{T}$ on $n$ vertices until $t$ time steps, then the geodesic distance $\mathcal{S}^{(w)}_{m}(t)$ of the $w$-star-fractal tree $\mathcal{T}^{\star}(t;w,m)$ will satisfy the following corollary.

\emph{Corollary 5} For a tree $\mathcal{T}$ on $n$ vertices, the exact solution for geodesic distance on its $w$-star-fractal tree $\mathcal{T}^{\star}(t;w,m)$ is

\begin{equation}\label{Section-3-c5}
\begin{aligned}\mathcal{S}^{(w)}_{m}(t)&=\Psi_{1}(m,w)\mathcal{S}^{(w)}_{m}(t-1)-\Psi_{2}(m,w)|\mathcal{V}^{(w)}_{m}(t-1)|^{2}+\Psi_{3}(m,n)|\mathcal{V}^{(w)}_{m}(t-1)|-\Psi_{4}(m,w)\\
&=[\Psi_{1}(m,w)]^{t}\mathcal{S}-\Psi_{4}(m,w)\sum_{i=0}^{t-1}[\Psi_{1}(m,w)]^{i}\\
&-\Psi_{2}(m,w)\sum_{i=0}^{t-1}[\Psi_{1}(m,w)]^{i}\left[(mw)^{t-1-i}\left(n-\frac{mw}{mw-1}\right)+\frac{mw}{mw-1}\right]^{2}\\
&+\Psi_{3}(m,w)\sum_{i=0}^{t-1}[\Psi_{1}(m,w)]^{i}\left[(mw)^{t-1-i}\left(n-\frac{mw}{mw-1}\right)+\frac{mw}{mw-1}\right].
\end{aligned}
\end{equation}

Here we omit the proof of Eq.(\ref{Section-3-c5}) with respect to the only requirement dependent of connection between  $|\mathcal{V}^{(w)}_{m}(t)|=(mw)^{t}\left(n-\frac{mw}{mw-1}\right)+\frac{mw}{mw-1}$ and Eq.(\ref{Section-3-t2}) in an iteration fashion.

\section{SEVERAL TREELIKE MODELS}

To show some applications of our technique to computations of geodesic distance over many treelike models, first, we will, in this section, make a brief description of generations of these such treelike models. In practice, some of them have been studied in more detail by virtue of many other methods, for instance, eigenvalue based on Laplacian matrix of underlying graph from spectral graph theory. As known, this is in principle a useful tool which can be widely used to capture a precise formula for geodesic distance over an arbitrary graph although sometimes some complicated calculations are involved. On the other hand, for the sake of simplicity, such a general method seems not perfect in comparison with some specific skills for exactly calculating geodesic distance on a number of graphs, for example, those with self-similarity. Referring to Ockham's razor, it is important and necessary to develop some new and professional methods for answering some special cases when people attempt to find out more general ways for addressing some kind of scientific issues. Indeed, our technique is built by analyzing the topological structure of treelike models in question thoroughly.

As we will point out shortly, some of the results below may be immediately viewed as consequences of the demonstrations in this preceding section. In other words, their correspondence underlying graphs are the special
cases of our formalism. In addition, while the other cases are obtained by no directly adopting our theorems, the lights used to accomplish the exact solution for them is indeed shed by our techniques. Put this another way, these techniques due to us are not limited to computations of geodesic distance over treelike models presented in this above section. Now let a well studied treelike model, the so-called T-graph, be as a starting point. We further develop the general foundations that will allow us
to handle the left two models, Cayley tree $C(t;n)$ and Exponential tree $\mathcal{T}(t;m)$, in stages.

\subsection{T-graph}

Technically speaking, the T-graph is established in an iterative manner. The seed is an edge connecting a pair of vertices, generally called $T_{0}$. As Fig.4(a) plots, one may apply $T-operation$ to $T_{0}$ and then obtain the next model $T_{1}$. Similarly, for time step $t\geq2$, the new generation $T_{t}$ can be obtained from its ancestor $T_{t-1}$ by implementing $T-$operation on each edge. In order to understand such a procedure, the rightmost of Fig.4 show an illustration of model $T_{3}$. It is easy to understand that $T-$operation can be naturally reduced as a special case of the ($1,m$)-star-fractal operation when parameter $m$ is assumed equal to $1$. Due to the fact, the ($1,m$)-star-fractal $\mathcal{T}^{\star}(t;1,m)$ can be viewed as a generalized version of T-graph where the seed is not only a single edge rather than an arbitrary tree. In principle, the seed may be a collection of trees, namely forest, as well.

\begin{figure}
\centering
\includegraphics[height=5cm]{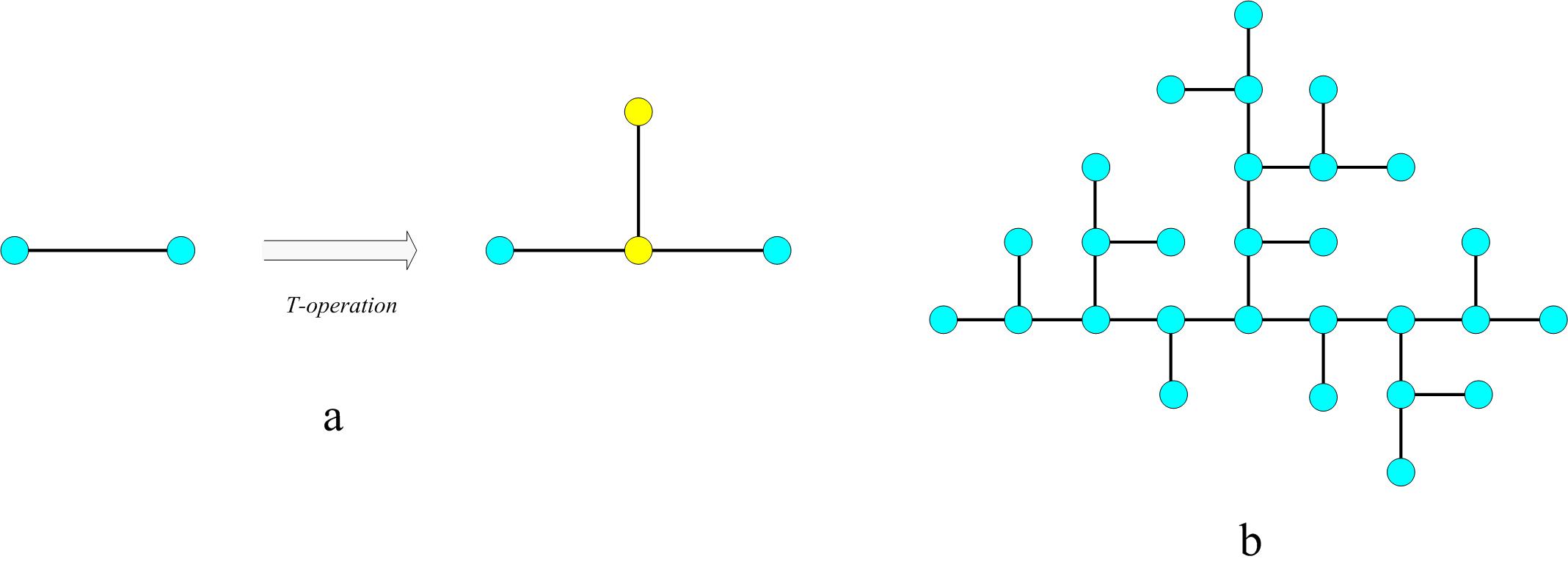}\\
{\small Fig.4. The diagram of \emph{T-operation} is plotted in panel (a). It is natural to see \emph{T-operation} is the most special case of our ($w,m$)-star-fractal operation where parameters $w$ and $m$ both equal to $1$. The rightmost is the diagram of T-graph at $t=3$.         }
\end{figure}

Obviously, after setting parameters $n=2$ and $m=1$ in Eq.(\ref{Section-3-c4}), the exact solution for geodesic distance $\mathcal{S}_{t}$ over T-graph $T_{t}$ will be written in the following way

\begin{equation}\label{Section-4-10}
\mathcal{S}_{t}=3^{t}+\frac{2^{t+2}+5}{5}\times3^{2t-1}-\frac{3^{t+1}}{5}.
\end{equation}

There is a long history of study on T-graph. One of significant reasons for this is that T-graph shows fractal property popularly seen in nature, for instance, Sierpinski gasket, Apollonian triangle, Koch curve, and so forth. Compared to these proposed fractal models, the T-graph is a fractal with the fractal dimension $d_{f}=\frac{\ln3}{\ln2}$. By using the same computation, we may know that the fractal dimensions of both ($1,m$)-star-fractal $\mathcal{T}^{\star}(t;1,m)$ and ($w,m$)-star-fractal $\mathcal{T}^{\star}(t;w,m)$ are, respectively, $d^{\star}_{f_{1,m}}=\frac{\ln(m+2)}{\ln2}$ and $d^{\star}_{f_{w,m}}=\frac{\ln[w(m+1)+1]}{\ln(w+1)}$. From the fractal dimension point of view, the both types of treelike modes have a choice to display the same fractal dimension, i.e.,

\begin{equation}\label{Section-4-11}
\frac{\ln(m+2)}{\ln2}=\frac{\ln[w(n+1)+1]}{\ln(w+1)}.
\end{equation}

A group of available values for Eq.(\ref{Section-4-11}) are $w=3,n=4$, and $m=2$, that is to say, $d^{\star}_{f_{1,2}}=d^{\star}_{f_{3,4}}=2$. In essence, there are a great deal of reasonable values satisfying Eq.(\ref{Section-4-11}). Interested reader can make some attempt to look for some intriguing solutions. Here we provide an equivalent condition of Eq.(\ref{Section-4-11}), i.e., $m=2^{\log_{(w+1)}[w(n+1)+1]}-2$. At the same time, Fig.5 can serve as an illustration of distribution of integer values with some constraints.

\begin{figure}
\centering
\includegraphics[height=5cm]{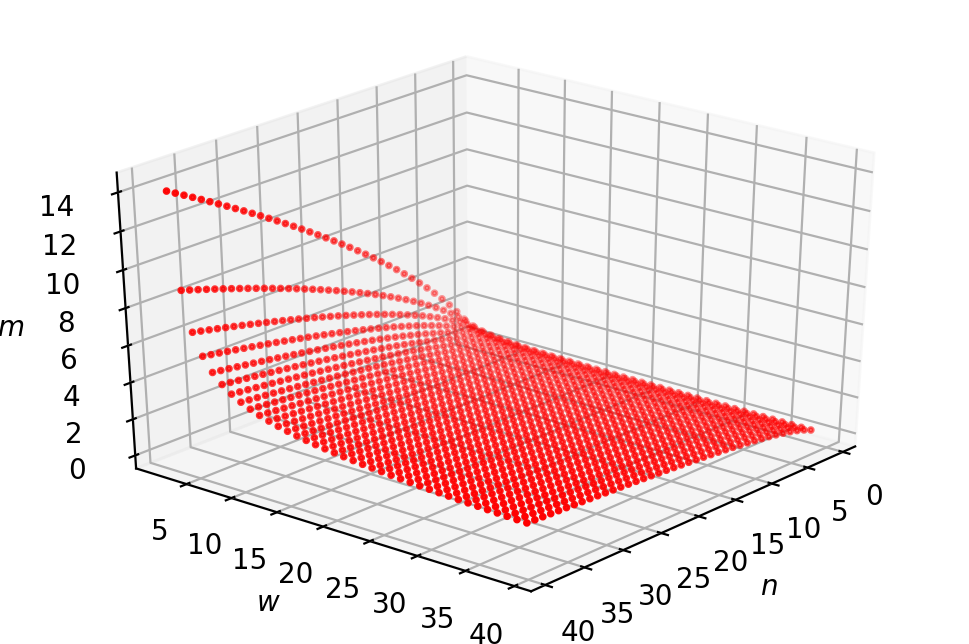}\\
{\small Fig.5. The diagram of solution for Eq.(\ref{Section-4-11}).         }
\end{figure}

By far, the operations discussed above are carried out on each edge of tree $\mathcal{T}$. By contrast, the next two treelike models will be generated in another growth manner, manipulation implemented on vertex of tree $\mathcal{T}$.

\subsection{Cayley tree $C(t;n)$}

 As we have seen, the models with fractal feature have self-similarity. Nonetheless, the converse of this statement is always not true. There is a great number of proofs verifying this point. Below is a better example, the so-called Cayley tree $C(t;n)$. First, to guarantee the self-contained nature of this paper, we will construct Cayley tree $C(t;n)$ by an algorithm.

\subsubsection{Construction of Cayley tree $C(t;n)$}

\emph{Algorithm 1}

At time step $t=1$, the seed $C(1;n)$ is a star with $n$ leaf vertices, as the star outlined in a yellow circle (online) in Fig.6 where there is in fact $4$ leaf vertices.

At time step $t=2$, the second tree $C(2;n)$ can be generated from model $C(1;n)$ by only connecting $n-1$ new vertices to each leaf vertex.

At time step $t\geq3$, the new generation tree $C(t;n)$ can be obtained from its ancestor model $C(t-1;n)$ using a similar manner to that mentioned above. A member of Cayley tree $C(t;n)$ is plotted in Fig.6.

Some of topological properties of Cayley tree $C(t;n)$ have been reported in numbers of literature. Here we just list part of their own rich properties. By using such a construction algorithm, it is clear to the eye that all nodes in the same shell are equivalent and have a unique degree. The nodes in the outermost shell have a degree equal to $1$, and all other nodes share a degree $n$ in common. Therefore, we can find that the vertex number of Cayley tree $C(t;n)$ is

\begin{equation}\label{Section-4-20}
|V_{C(t;n)}|=\frac{n(n-1)^{t}-2}{n-2}.
\end{equation}

It is worthy noting that although Cayley tree $C(t;n)$ is obviously self-similar, its fractal
dimension is infinite, and so it has no fractal character. Just as with discussions about geodesic distance over some graphs having self-similarity, we will in the next subsection focus on computation of geodesic distance over Cayley tree $C(t;n)$ in a recursive fashion. Some additional notations will be taken to help us completely accomplish the whole computation and their concrete meanings will be given at the positions of their appearance.

\begin{figure}
\centering
\includegraphics[height=5cm]{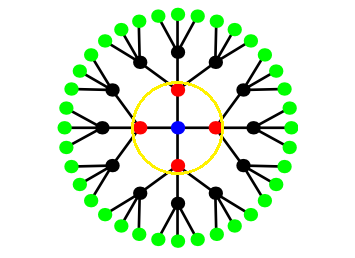}\\
{\small Fig.6. The diagram of treelike model $C(t,n)$ that is commonly called Cayley tree where $n=4$ and $t=3$. For our purpose, we would like to employ another description for Cayley tree $C(t,n)$ in which the original star made up one blue vertex and $n$ red vertices is defined as the seed of our new representation $C^{\ast}(t,\psi)$ where $\psi=n-1=3$.  }
\end{figure}

\subsubsection{Calculation of geodesic distance upon self-similarity}

First of all, we would like to make an abstract description of Cayley tree $C(t,n)$ in another manner due to its self-similarity. Cayley tree $C(t,n)$ can be reconstructed by merging the root of each of $n$ smaller ingredients $A_{t}$ into a key vertex as shown in Fig.7. Each ingredient is a rooted tree whose root is the central vertex of the original star. Based on such a construction, we have

\begin{equation}\label{Section-4-21}
|V_{C(t,n)}|=n|A_{t}|-(n-1)
\end{equation}
where notation $|A_{t}|$ represents the vertex number of ingredient $A_{t}$. By the same trick, the expression of geodesic distance over Cayley tree $C(t,n)$ can be certainly written as

\begin{equation}\label{Section-4-22}
\mathcal{S}_{C(t,n)}=n\Gamma_{A_{t}}+\Omega_{C(t,n)}.
\end{equation}

The first term $\Gamma_{A_{t}}$ of the right-hand side of Eq.(\ref{Section-4-22}) is the solution for geodesic distance over ingredient $A_{t}$. The other term is the sum of geodesic distance over all possible pairs of ingredients and thus can easily be expressed

\begin{equation}\label{Section-4-23}
\Omega_{C(t,n)}=\sum_{1\leq i<j\leq n}\Omega^{ij}_{C(t,n)}=\frac{n(n-1)}{2}\Omega^{12}_{C(t,n)}
\end{equation}
where we have used the fact that ingredients are self-similar with each other.

By definition, the sub-formula of Eq.(\ref{Section-4-23}) is able to be translated as follows

\begin{equation}\label{Section-4-24}
\begin{aligned}\Omega^{12}_{C(t,n)}&=\sum_{i\in A^{1}_{t},i\neq \omega^{1}_{t}}\sum_{j\in A^{2}_{t},j\neq \omega^{2}_{t}}d^{ij}_{t}=\sum_{i\in A^{1}_{t},i\neq \omega^{1}_{t}}\sum_{j\in A^{2}_{t},j\neq \omega^{2}_{t}}(d^{i\theta_{t}}_{t}+d^{\theta_{t}j}_{t})=2(|A_{t}|-1)\Theta_{C(t,n)}
\end{aligned}
\end{equation}
in which symbol $\omega^{i}_{t}$ is the root of ingredient $ A^{i}_{t}$ and becomes in essence the central vertex of Cayley tree $C(t,n)$ after merging operation. The symbol $\theta_{t}$ is, thus, the central vertex. Furthermore, the sum of geodesic distances over all possible pairs of vertices where one of which is vertex $\omega^{i}_{t}$ and the other is an arbitrary vertex belonging to $ A^{i}_{t}$ which is distinct with vertex $\omega^{i}_{t}$.

Armed with self-similar feature and Eqs.(\ref{Section-4-23})-(\ref{Section-4-24}), the $\Gamma_{A_{t}}$ of Eq.(\ref{Section-4-22}) may be reorganized in the following way

\begin{figure}
\centering
  \includegraphics[height=5cm]{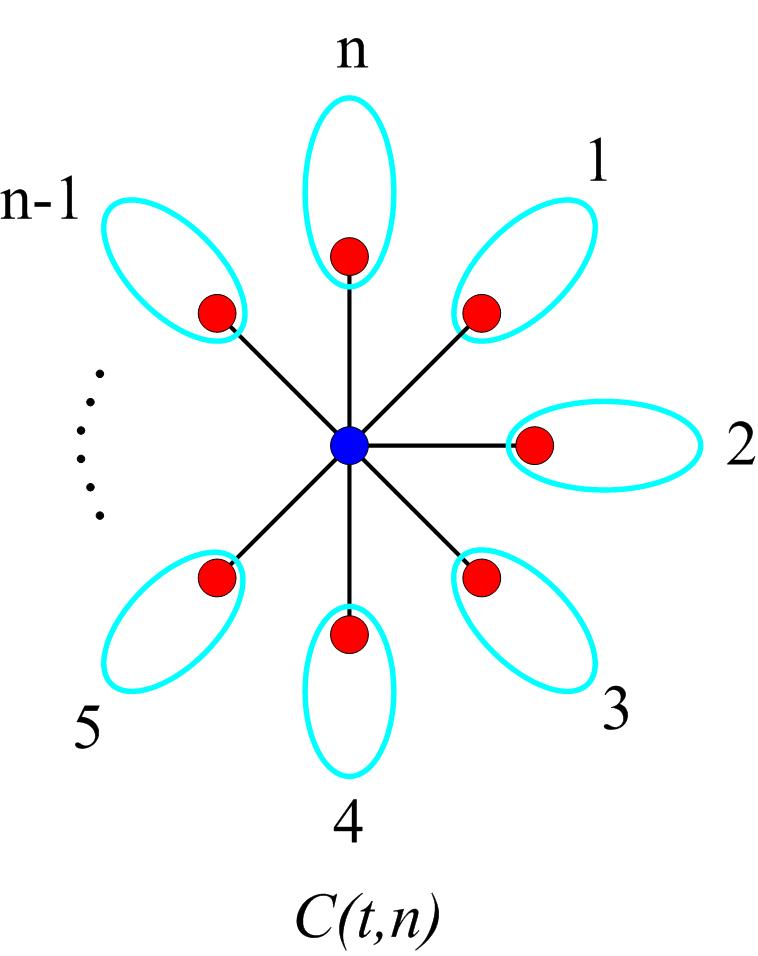}\\
{\small Fig.7. The diagram of Cayley tree $C(t,n)$. We denote by the collection of the blue vertex attached to an indigo ellipse by an edge a smaller ingredient $A_{t}$ which is in fact a rooted tree.  }
\end{figure}

\begin{equation}\label{Section-4-25}
\Gamma_{C(t,n)}=(n-1)\Gamma_{C(t-1,n)}+\frac{(n-1)(n-2)}{2}\Omega^{12}_{C(t-1,n)}+\Theta_{C(t,n)}.
\end{equation}

Similarly, we can find

\begin{equation}\label{Section-4-26}
\Theta_{C(t,n)}=(n-1)\Theta_{C(t-1,n)}+|A_{t}|-1.
\end{equation}

Now, the determination of Eq.(\ref{Section-4-26}) depends on determining the value of $|A_{t}|$. Solving for $|A_{t}|$ through Eq.(\ref{Section-4-20}) and Eq.(\ref{Section-4-21}) yields

 \begin{equation}\label{Section-4-27}
|A_{t}|=\frac{(n-1)^{t}+n-3}{n-2}.
\end{equation}

Substituting Eq.(\ref{Section-4-27}) into Eq.(\ref{Section-4-26}) and making use of necessary arithmetics together produces

 \begin{equation}\label{Section-4-28}
\Theta_{C(t,n)}=(n-1)^{t-1}\Theta_{C(1,n)}+\sum_{i=0}^{t-2}(n-1)^{i}\times\frac{(n-1)^{t-i}-1}{n-2}.
\end{equation}

With the initial condition $\Theta_{C(1,n)}=1$, the value of $\Theta_{C(t,n)}$ is

 \begin{equation}\label{Section-4-29}
\Theta_{C(t,n)}=\frac{[(n-2)t-1](n-1)^{t}+1}{(n-2)^{2}}.
\end{equation}

Plugging both Eq.(\ref{Section-4-27}) and Eq.(\ref{Section-4-29}) into Eq.(\ref{Section-4-24}) outputs

 \begin{equation}\label{Section-4-210}
\Omega^{12}_{C(t,n)}=\frac{2}{(n-2)^{3}}([(n-2)t-1](n-1)^{2t}-[(n-2)t-2](n-1)^{t}-1).
\end{equation}

By analogy with the forgoing calculations,  Eq.(\ref{Section-4-29}) and Eq.(\ref{Section-4-210}) will be utilized to determine the closed-form expression of $\Gamma_{C(t,n)}$ in Eq.(\ref{Section-4-25}) as below

 \begin{equation}\label{Section-4-211}
\Gamma_{C(t,n)}=(n-1)^{t-1}+\frac{(t-1)(n-1)^{t+1}-(n-1)^{t-1}+1}{(n-2)^{2}}+\frac{[(n-2)t-n](n-1)^{2t}+2(n-1)^{t+1}}{(n-2)^{3}}.
\end{equation}

Finally, with the help of Eq.(\ref{Section-4-23}), Eq.(\ref{Section-4-210}) and Eq.(\ref{Section-4-211}), the exact solution for geodesic distance over Cayley tree $C(t,n)$ as shown in Eq.(\ref{Section-4-22}) can be obtained

 \begin{equation}\label{Section-4-212}
\begin{aligned}\mathcal{S}_{C(t,n)}&=n\Gamma_{C(t,n)}+\Omega_{C(t,n)}\\
&=n(n-1)^{t-1}+\frac{n(t-1)(n-1)^{t+1}-n(n-1)^{t-1}+n}{(n-2)^{2}}\\
&+\frac{n[(n-2)t-n](n-1)^{2t}+2n(n-1)^{t+1}}{(n-2)^{3}}\\
&+\frac{n(n-1)}{(n-2)^{3}}([(n-2)t-1](n-1)^{2t}-[(n-2)t-2](n-1)^{t}-1).
\end{aligned}
\end{equation}

We remark that some evidently conditional values have been carefully considered here, for instance, $|A_{1}|=2, |A_{2}|=n+1$, $\Theta_{C(1,n)}=1$ and $\Omega^{12}_{C(1,n)}=2$, as well $\Gamma_{C(1,n)}=1$.

We must stress again that the self-similarity of Cayley tree $C(t,n)$ plays a vital role on successfully finishing the computation of geodesic distance. Nevertheless, it is clear to see that the seed of Cayley tree $C(t,n)$ is a star with $n$ leaf vertices. Interesting although this type of Cayley tree $C(t,n)$ is, however, it is a very special case. If one would to like to generalize the formulism, then self-similarity remains unchanged in the limit of large graph size but the method presented above assumes not enough convenient to utilize. A simple choice is to suppose that the seed is an arbitrary tree with some constraints . To satisfy the requirement of the intrinsic property of Cayley tree $C(t,n)$, i.e., all but leaf vertices having a unique degree equal to $n$, a candidate may be easily generated by applying ($w,m$)-star-fractal operation on a single edge where the only thing is to set $m=n-2$. More generally, it is essential to develop any other tools for addressing these such problems. Surprisingly, the lights shed by the building-up of the above two theorems will act as a guide to help us achieve our desired goal.

Last but not least, let us return our insights to the present questions. It should be worthy mentioning that we will take useful advantage of another description of the generalized version of Cayley tree $\mathbf{C}(t,n)$, as illustrated in Fig.(6).

\subsubsection{Another method for calculation of geodesic distance}

Consider an arbitrary tree $\mathcal{T}$ on $|V|$ with requirement as the seed of Cayley tree $\mathbf{C}(t,n)$, we may suppose that the geodesic distance over tree $\mathcal{T}$ is known and assumed equivalent to $\mathcal{S}$. As before, for the new generation tree $\mathbf{C}(1,n)$, its corresponding geodesic distance $\mathcal{S}(1)$ should follow

 \begin{equation}\label{Section-4-220}
\left\{\begin{aligned}&\mathcal{S}(1)=\mathcal{S}+S_{1}(1)+S_{2}(1)\\
&S_{1}(1)=\psi S_{1}(0)+2\psi|\Delta V_{0}|+\psi S_{2}(0)\\
&S_{2}(1)=\psi ^{2}S_{2}(0)+\psi ^{2}|\Delta V_{0}|^{2}-\psi|\Delta V_{0}|
\end{aligned}
\right.
\end{equation}
here parameter $\psi=n-1$ is used and we denote by $\Delta V_{t_{i}}$ the leaf vertex set of Cayley tree $\mathbf{C}(t_{i},n)$. Meantime, let $S_{2}(t_{i})$ represent the sum of geodesic distance over all pairs of leaf vertices in $\Delta V_{t_{i}}$ and $S_{1}(t_{i})$ be equal to the sum of geodesic distance over all pairs of vertices one of which is from $\Delta V_{t_{i}}$ and the other from $V_{\mathbf{C}(t_{i},n)}-\Delta V_{t_{i}}=V_{\mathbf{C}(t_{i}-1,n)}$ ($t_{i}\geq1$).

After $t$ time steps, the geodesic distance $\mathcal{S}(t)$ over the generalized Cayley tree $\mathbf{C}(t,n)$ can be expressed in an iteration fashion
\begin{equation}\label{Section-4-221}
\begin{aligned}
\mathcal{S}(t)&=\mathcal{S}+\sum_{i=1}^{2}\sum_{j=1}^{t}S_{i}(j)=\mathcal{S}+\sum_{i=1}^{t}\left(\psi^{i}S_{1}(0)+\sum_{j=0}^{i-1}2\psi^{i-j}|\Delta V_{j}|\right)\\
&+\sum_{i=1}^{t}\left(\psi^{2i}S_{2}(0)+\sum_{j=0}^{i-1}\psi^{2j+2}|\Delta V_{i-1-j}|^{2}-\sum_{j=0}^{i-1}\psi^{2j+1}|\Delta V_{i-1-j}|\right)\\
&+\sum_{l=1}^{t}\sum_{i=0}^{l-1}\psi^{l-i}\left(\psi^{2i}S_{2}(0)+\sum_{j=0}^{i-1}\psi^{2j+2}|\Delta V_{i-1-j}|^{2}-\sum_{j=0}^{i-1}\psi^{2j+1}|\Delta V_{i-1-j}|\right).
\end{aligned}
\end{equation}

Apparently, the technique described here is more general and hence enough light to adopt to exactly obtain the closed-form solution for geodesic distance over the generalized Cayley tree $\mathbf{C}(t,n)$, including that special case of a single edge chosen as a seed. The only thing to really note here is to know the initial value $\mathcal{S}$ for geodesic distance over the seed in advance before starting with the latter part of computation. Compared to the vertex number $|V_{\mathbf{C}(t,n)}|$ of the generalized Cayley tree $\mathbf{C}(t,n)$ increasing exponentially as a function of time step $t$, it is not hard to derive a precise expression of $\mathcal{S}$ on a relatively smaller seed using exhaustively enumeration method on current computer.

As an obvious consequence of application of Eq.(\ref{Section-4-221}), suppose that $|\Delta V_{1}|=\psi+1$, $\mathcal{S}(1)=(\psi+1)^{2}$, $S_{1}(1)=\psi+1$ and $S_{2}(1)=\psi(\psi+1)$, one can gain the same result as said by Eq.(\ref{Section-4-212}). In practice, on the basis of the above assumption, we generates the Cayley tree $C(t,n)$ as output by Algorithm 1.

\subsection{Exponential tree $\mathcal{T}(t;m)$}

Till now, all the treelike models share a feature in common, i.e., all except for leaf vertices with a unchanged degree. To distinguish the difference between these models and the latter treelike models as we will introduce shortly, the former is considered homogeneous and the latter is heterogeneous. Since then, we will propose a class of treelike models, hereafter called Exponential tree $\mathcal{T}(t;m)$, which have been studied in our previous work \cite{X-M-2019}. Different from the published results, here our description may be thought of as a more general form mainly because the seed is an arbitrary tree $\mathcal{T}$ we please.

\emph{Algorithm 2}

At $t=0$, the seed $\mathcal{T}$ is a tree on $|V|$ vertices, which is also regarded as $\mathcal{T}(0,m)$.

At $t=1$, the next model $\mathcal{T}(1,m)$ can be generated by connecting $m$ new vertices to each vertex model $\mathcal{T}(0,m)$.

At $t\geq2$, the second model $\mathcal{T}(t,m)$ can be constructed in a similar manner to that mentioned at the preceding time step.

After $t$ time steps, the vertex number of Exponential tree $\mathcal{T}(t,m)$ obeys
\begin{equation}\label{Section-4-31}
|V(t,m)|=(m+1)^{t}|V|.
\end{equation}

By analogy with the development of Eqs.(\ref{Section-4-220}) and (\ref{Section-4-221}), the geodesic distance $\mathcal{S}(t,m)$ over Exponential tree $\mathcal{T}(t,m)$ follows
\begin{equation}\label{Section-4-32}
\begin{aligned}
\mathcal{S}(t,m)&=(1+m)^{2t}\mathcal{S}+m(m+1)\sum_{i=0}^{t-1}(1+m)^{2i}|V(t-i-1,m)|^{2}-m\sum_{i=0}^{t-1}(1+m)^{2i}|V(t-i-1,m)|
\end{aligned}
\end{equation}
in which $\mathcal{S}$ represents the known geodesic distance over the seed $\mathcal{T}$.

We remark that the reason why we call the new generation model $\mathcal{T}(t,m)$ Exponential tree is that its degree distribution takes an exponential form in the limit of large graph size, i.e., $P_{cum}(k\geq k_{t_{i}})\sim \text{exp}[-(m+1)k_{t_{i}}]$.

As with discussions on complex network study, let the seed be a single edge, then using Eq.(\ref{Section-4-32}) we can obtain a consequence

\begin{equation}\label{Section-4-33}
\mathcal{S}_{t}=(m+1)^{t-1}[2+(4mt+m-1)(m+1)^{t}]
\end{equation}
as reported in \cite{Z-2010}.

\section{RANDOM WALKS ON TREELIKE MODELS}

As the discrete-time representative of Brownian motion and diffusive processes, random walk, which is put first forward by Einstein and Smoluchowski, have proven useful in a wide range of applications. Particularly, the recent two decades witnesses an upsurge of understanding random walk on complex networked models particularly because of the hyperactive complex network study \cite{Alfred-2018}-\cite{ Rushabh-2016}. Without loss of generality, we still denote by $G(V,E)$ a complex networked model.

Random walk in fact describes an ideal situation in which a walker (particle) has a lack of any information of the underlying graph $G(V,E)$ and just chooses uniformly at random one vertex of its neighbor set to move on, called the unbiased Markov random walks as well. The one of most significant issues of studying random walk on a given model $G(V,E)$ is determining precise solutions of mean first-passage time ($MFPT$). This quantity is by definition equivalent to the average over $FPT$s of all possible pairs of vertices $u$ and $v$ in $V$. Generally speaking, the $FPT$ for any vertex pair can be determined on the basis of the fundamental matrix corresponding to the underlying graph $G(V,E)$. Facing with a graph with thousands of vertices, the matrix-based method will be prohibitively complicated. However, for a graph with tree form, the above task remains solvable in view of a lemma built upon effective resistance of electrical network as we will state below.

\textbf{Lemma 3} For a given tree $\mathcal{T}$ on $|V|$ vertices, there is a close relationship between the geodesic distance $\mathcal{S}$ over it and the mean first-passage time $MFPT$ for a random walker on the tree

\begin{equation}\label{Section-5-1}
MFPT=\frac{\mathcal{S}}{|V|}.
\end{equation}

For the sake of space limitation, we omit the proof of Eq.(\ref{Section-5-1}) and interested reader is referred to \cite{M-W-L-2019} for more detail. In other words, determining the solution for geodesic distance on a tree $\mathcal{T}$also determines the mean first-passage time $MFPT$ for a random walker on this tree. Keep this in mind, we will give some obvious propositions by virtue of Eq.(\ref{Section-5-1}). It is worthy noting that the below results hold in the limit of large graph size. Hence, in such a situation, the seed may be a single edge.

\textbf{Proposition 1} For the $m$th-order subdivision tree $\mathcal{T}^{m}(t)$, its mean first-passage time $MFPT_{\mathcal{T}^{m}(t)}$ is

\begin{equation}\label{Section-5-2}
MFPT_{\mathcal{T}^{m}(t)}\sim |V_{\mathcal{T}^{m}(t)}|^{2}.
\end{equation}

\textbf{Proposition 2} For the ($1.m$)-star-fractal tree $\mathcal{T}^{\star}(t;1,m)$, its mean first-passage time $MFPT_{\mathcal{T}^{\star}(t;1,m)}$ is

\begin{equation}\label{Section-5-2}
MFPT_{\mathcal{T}^{\star}(t;1,m)}\sim |V_{\mathcal{T}^{\star}(t;1,m)}|^{\frac{\ln2(m+2)}{\ln2}}.
\end{equation}

\textbf{Proposition 3} For the ($w.m$)-star-fractal tree $\mathcal{T}^{\star}(t;w,m)$, its mean first-passage time $MFPT_{\mathcal{T}^{\star}(t;w,m)}$ is

\begin{equation}\label{Section-5-2}
MFPT_{\mathcal{T}^{\star}(t;w,m)}\sim |V_{\mathcal{T}^{\star}(t;w,m)}|^{\frac{\ln(w+1)[m(w+1)+1]}{\ln(w+1)}}.
\end{equation}

\textbf{Proposition 4} For the T-graph $T_{t}$, its mean first-passage time $MFPT_{T_{t}}$ is

\begin{equation}\label{Section-5-2}
MFPT_{T_{t}}\sim |V_{T_{t}}|^{\frac{\ln6}{\ln2}}.
\end{equation}

\textbf{Proposition 5} For the Cayley tree $C(t,n)$, its mean first-passage time $MFPT_{C(t,n)}$ is

\begin{equation}\label{Section-5-2}
MFPT_{C(t,n)}\sim \Upsilon_{1}(n,t)|V_{C(t,n)}|,\quad\quad \Upsilon_{1}(n,t)= O\left(\left(\frac{n}{n-2}\right)^{2}t\right).
\end{equation}

\textbf{Proposition 6} For the Exponential tree $\mathcal{T}(t;m)$, its mean first-passage time $MFPT_{\mathcal{T}(t;m)}$ is

\begin{equation}\label{Section-5-2}
MFPT_{\mathcal{T}(t;m)}\sim \Upsilon_{2}(m,t)|V_{\mathcal{T}(t;m)}|,\quad\quad \Upsilon_{2}(m,t)= O\left(\frac{4m}{m+1}t\right).
\end{equation}

The propositions obtained above show our desired results. The heuristic explanations about how the operations discussed in this paper have a considerable influence on $MFPT$ over underlying graphs will be studied in depth in the following section.

\section{DISCUSSION}

This section is intended as an exploration to similarity and difference among the operations introduced here. To this end, the available classifications is on demand. In so doing, the below discussions will work well in stages. As before, we remak that all prerequisites will be initialized, i.e., 1) a single edge chosen as the seed, 2) the large value for $t$.

\textbf{Statement 1 } All the generated treelike models have obviously self-similarity.

\textbf{Statement 2 } Here, we need to take advantage of a notation, denoted by $\Delta V$, which suggests how fast the ratio of vertex number change is. Technically speaking, $\Delta V$ is defined as the ratio $\lim_{t\rightarrow\infty}(|V_{t}|-|V_{t-1}|)/|V_{t-1}|$. As before, we list some results in order.
$$\Delta V_{\mathcal{T}^{m_{1}}(t)}=m_{1},\quad \Delta V_{\mathcal{T}^{\star}(t;1,m_{2})}=m_{2}+1,\quad \Delta V_{\mathcal{T}^{\star}(t;w,m_{3})}=w(m_{3}+1),$$
$$\Delta V_{T_{t}}=3,\quad \Delta V_{C(t,n)}=n-2,\quad \Delta V_{\mathcal{T}(t;m_{4})}=m_{4}.$$

In some cases, the five of which can have a common $\Delta V$. Nonetheless, it is easy to see that the five resulting treelike models are different from each other only because of completely distinct topology structure. This implies that a unique $\Delta V$ is not necessary to result in a unique structure. More generally, this leads to a significant difference from the topological structure point of view as we will state shortly.

\textbf{Statement 3 } Different from statement 1, graph with fractal property should be self-similar. On the other hand, graph with self-similarity is not necessarily fractal. Here, three treelike models, i.e., the $m$th-order subdivision tree $\mathcal{T}^{m}(t)$, the Cayley tree $C(t,n)$ and the Exponential tree $\mathcal{T}(t;m)$, all have no fractal property. The rest three ones, the ($1.m$)-star-fractal tree $\mathcal{T}^{\star}(t;1,m)$, the ($w.m$)-star-fractal tree $\mathcal{T}^{\star}(t;w,m)$ together with the T-graph $T_{t}$, are all fractal.

Meantime, for the T-graph $T_{t}$, its random-walk dimension $d_{T_{t}}$ satisfies $d_{w_{T_{t}}}=\frac{\ln6}{\ln2}=1+d_{f_{T_{t}}}$ where $d_{f_{T_{t}}}$ represents the fractal dimension of T-graph $T_{t}$. Similarly, the ($1.m$)-star-fractal tree $\mathcal{T}^{\star}(t;1,m)$ and the ($w.m$)-star-fractal tree $\mathcal{T}^{\star}(t;w,m)$ have random-walk dimension $d_{w_{\mathcal{T}^{\star}(t;1,m)}}=\frac{\ln2(m+2)}{\ln2}=1+d_{f_{\mathcal{T}^{\star}(t;1,m)}}$ and $d_{w_{\mathcal{T}^{\star}(t;w,m)}}=\frac{\ln(w+1)[m(w+1)+1]}{\ln(w+1)}=1+d_{f_{\mathcal{T}^{\star}(t;w,m)}}$, respectively.

Referring to statement in \cite{Erik-2005}, one can estimate whether the random walk on a graph with fractal feature is persistent or not with respect to the spectral dimension $\widetilde{d}=\frac{2d_{f}}{d_{w}}$. In another word, the random walk phenomenon on such a graph is termed ``persistence'' when $\widetilde{d}<2$ is true, is not otherwise. Using the results above, it is easy to compute

$$\widetilde{d_{T_{t}}}=\frac{2d_{f_{T_{t}}}}{d_{w_{T_{t}}}}=\frac{\ln9}{\ln6}<2, \quad\widetilde{d_{\mathcal{T}^{\star}(t;1,m)}}=\frac{2d_{f_{\mathcal{T}^{\star}(t;1,m)}}}{d_{w_{\mathcal{T}^{\star}(t;1,m)}}}=\frac{\ln(m+2)^{2}}{\ln2(m+2)}<2,$$

$$\widetilde{d_{\mathcal{T}^{\star}(t;w,m)}}=\frac{2d_{f_{\mathcal{T}^{\star}(t;w,m)}}}{d_{w_{\mathcal{T}^{\star}(t;w,m)}}}=\frac{\ln[m(w+1)+1]^{2}}{\ln(w+1)[m(w+1)+1]}<2.$$

Obviously, a walker starting from a given node of the three treelike models above will return to the node almost surely over the course of time. In addition, the left three have no such properties.

From the scaling point of view, the $MFPT$ of the former four treelike models are all a power-law function of its vertex number, separately. On the other hand, the latter two are not power-law but a compound product of two functions, i.e., one for time step $t$ and the other for vertex number. This further implies that topological structure has a significant impact on the $MFPT$ for  random walk on its correspondence underlying graph.

\textbf{Statement 4 } As another metric, the mean geodesic distance can be utilized to determine how fast information diffusion on an underlying network. By Eq.(\ref{Section-11}), we immediately give a series of consequences in the large graph size limit
$$\langle\mathcal{S}_{\mathcal{T}^{m}(t)}\rangle\sim |V_{\mathcal{T}^{m}(t)}|, \quad \langle\mathcal{S}_{\mathcal{T}^{\star}(t;1,m)}\rangle\sim |V_{\mathcal{T}^{\star}(t;1,m)}|^{d_{f_{\mathcal{T}^{\star}(t;1,m)}}},$$
$$ \langle\mathcal{S}_{\mathcal{T}^{\star}(t;w,m)}\rangle\sim |V_{\mathcal{T}^{\star}(t;w,m)}|^{d_{f_{\mathcal{T}^{\star}(t;w,m)}}}, \quad \langle\mathcal{S}_{T_{t}}\rangle\sim |V_{\mathcal{T}^{m}(t)}|^{d_{f_{T_{t}}}},$$
$$ \langle\mathcal{S}_{C(t,n)}\rangle\sim \ln|V_{C(t,n)}|, \quad \langle\mathcal{S}_{\mathcal{T}(t;m)}\rangle\sim \ln|V_{\mathcal{T}(t;m)}|.$$

\section{CONCLUSION}

To conclude, we first introduce some operations widely studied in a great variety of disciplines including the $m$th-order subdivision ($m\geq1$) and the ($w,m$)-star-fractal operation ($w\geq1,m\geq1$). These operations can be implemented on any graph. In particular, one of the most fundamental graphs, tree, is chosen as our focus. Using each of the operations above, we may obtain a resulting tree. On the basis of the resulting treelike models, we discuss an interesting and significant topological parameter, i.e., geodesic distance. This quantity is closely correlated with diffusion on underlying graph and hence may usually be an indicator of how fast information propagation is. To analytically compute the exact solution for geodesic distance on each treelike models, we do not use some more general methods, for instance, that matrix-based ones, rather than make use of mapping-based technique. And then, we indeed obtain our desired formulas. At the same time, our novel technique is more convenient to manipulate than that matrix-based ones, as least on the treelike models under consideration.

To further highlight convenience and usefulness of our technique, two types of treelike models of significant interest, Cayley tree $C(t,n)$ and Exponential tree $\mathcal{T}(t;m)$, are considered. Different from the previous discussions in some literature, the seed of the both models is not a single edge but a tree with some requirements. Put this another way, our generated models can be regarded as generalized versions. So, the pre-existing results is a special case of ours.

Finally, by using a lemma connecting the geodesic distance to mean first-passage time on tree, we give some intriguing propositions and then make a few heuristic explanations. Among of which different operations can have a vital influence on the two quantities as explained in section 6.

We would like to stress that our work is only a tip of the iceberg and however the lights shed by our
methods can be beneficial. We wish to witness some other applications of our methods in the days to come. Meantime, there still are some open questions which are waiting for us to address.

\section*{ACKNOWLEDGMENT}

The authors would like to thank Xudong Luo for useful conversations. The research was supported by the National Key Research and Development Plan under grant 2017YFB1200704 and the National Natural Science Foundation of China under grant No. 61662066.

\end{document}